\documentclass[a4paper,12pt]{article}

\usepackage{amsmath,amssymb}
\usepackage{graphicx}
\usepackage[a4paper,margin=2.5cm]{geometry}
\usepackage{tikz}
\usepackage{cite}
\usepackage[colorlinks=true,linkcolor=blue,citecolor=blue]{hyperref}

\newtheorem{thm}{{\bf T}{\footnotesize \bf HEOREM}}
\newtheorem{lm}[thm]{{\bf L}{\footnotesize \bf EMMA}}

\newtheorem{conj}[thm]{{\bf C}{\footnotesize \bf ONJECTURE}}

\newenvironment{Proof}
{\noindent {\it Proof.} }
{\hbox{\rule{6pt}{6pt}}\bigskip}

\newenvironment{Proofof}[1]
{\noindent {\it Proof of {#1}.} }
{\hbox{\rule{6pt}{6pt}} \bigskip}

\begin{document}
\sloppy
\title{Another proof of the result on rotation compatible planar covers}

\author{
Shohei Koizumi\thanks{%
Graduate School of Science and Technology,
Niigata University,
8050 Ikarashi 2-no-cho, Nishi-ku, Niigata, 950-2181, Japan.
E-mail: {\tt s-koizumi@m.sc.niigata-u.ac.jp}}, 
Yusuke Suzuki\thanks{
Department of Mathematics, Niigata University, 
8050 Ikarashi 2-no-cho, Nishi-ku, Niigata, 950-2181, Japan.
Email: {\tt y-suzuki@math.sc.niigata-u.ac.jp}} \ 
and Kensuke Tamura\thanks{
Graduate School of Science and Technology, 
Niigata University, 
8050 Ikarashi 2-no-cho, Nishi-ku, Niigata, 950-2181, Japan.
}
}

\date{}
\maketitle

\begin{abstract}
Negami's Planar Cover Conjecture asserts that a connected graph has a finite planar cover if and only if it 
can be embedded on the projective plane. 
While this statement has already been proven for rotation compatible planar covers, 
namely covers equipped with a certain condition on the rotation system, 
the existing proof relies on advanced algebraic and topological methods.
In this paper, we provide another proof of this result, 
focusing primarily on combinatorial arguments based on a structural analysis 
with respect to a spanning tree in the base graph. 
\end{abstract}


\section{Introduction}\label{intro}
In this paper, we consider only finite simple graphs.  
We denote the vertex set and edge set of a graph $G$ by 
$V(G)$ and $E(G)$, respectively. 
For a vertex $v$ of a graph $G$, 
the {\em neighborhood} of $v$ is denoted by $N_G(v)$, 
or simply $N(v)$ if the graph is clear from the context. 
For a $S\subset V(G)$, we denote the subgraph induced by $S$ 
by $G[S]$. 
For terminology not described here, we 
follow the standard graph theory terminology (e.g., see \cite{Ds}).
Let $\tilde{G}$ and $G$ be graphs. 
If there exists an ($n$-to-one) surjection $p \colon V(\tilde{G}) \to V(G)$ such that  
for any vertex $\tilde{v} \in V(\tilde{G})$ the restriction 
$p|_{N(\tilde{v})} \colon N_{\tilde{G}}(\tilde{v}) \to N_{G}(p(\tilde{v}))$ is a bijection, 
then $\tilde{G}$ is an ($n$-{\em fold\/}) {\em cover} of $G$. 
Then, $p$ is called a {\em covering map} (or {\em cover projection}) of the cover, 
and we also write $p \colon \tilde{G} \to G$ for such a covering map.
Let $L = (v_1, \ldots, v_k)$ be a list (or a cyclic order) of vertices of $\tilde{G}$. 
We define $p(L)$ to be the list (or cyclic order) $(p(v_1), \ldots, p(v_k))$; 
note that vertices may appear more than once in this list. 
The following is the well-known {\em Planar Cover Conjecture}, 
proposed in 1986.
\begin{conj}[Negami \cite{Negami1988PCC}]\label{conj:Negami}
A connected graph has a finite planar cover if and only if it can be embedded on the projective plane.
\end{conj} 

In the conjecture above, the sufficiency is clear, but the necessity is still open.
Unlike the conventional scheme which primarily focuses on showing that $K_{1,2,2,2}$ 
has no finite planar cover (see e.g., \cite{n14, Arch, Arch2, K44e, Hli, HandT, 20years, kgraphs}), 
in 2024, Negami published a paper \cite{Negami2024main} presenting a new approach 
to Conjecture~\ref{conj:Negami}, which focuses 
on the ``rotation system'' on vertices as follows. 

Let $G$ be a graph embedded on an orientable closed surface $F^2$, 
and we fix the clockwise orientation of $F^2$. 
For each vertex $v$ of $G$, the clockwise orientation induces 
a cyclic order on $N(v)$. 
We call the cyclic order above a {\em rotation} of $v$, and denote 
it by $\sigma_v$. 
Then $\sigma_v$ can be regarded as a map 
$\sigma_v \colon N(v) \to N(v)$ 
that sends each vertex $u \in N(v)$ to its successor 
around $v$ with respect to the clockwise orientation. 
Then $\sigma_v^{-1}$ is the cyclic order on $N(v)$ defined by the 
anticlockwise orientation around $v$; 
that is, the inverse of $\sigma_v$. 
The {\em rotation system} $\sigma$ is defined as 
$\sigma = \{ \sigma_v \colon v \in V(G)\}$. 

Let $G$ be a connected graph and let $\tilde G$ be a 
cover of $G$ with a covering map $p \colon \tilde{G} \to G$. 
Assume that $\tilde{G}$ is embedded on an orientable closed surface $F^2$, 
and has a rotation system $\sigma$ according to the clockwise orientation of $F^2$. 
Negami \cite{Negami2024main} proposed the following two assumptions given to 
the composite structure $(\tilde{G}, \sigma)$. 
In fact, $(\tilde{G}, \sigma)$ is {\em rotation compatible} if 
$(\tilde{G}, \sigma)$ satisfies those two assumptions.


\bigskip

\noindent
\textbf{Assumption 1:} 
For any two vertices $v_1$ and $v_2$ of $\tilde{G}$ such that $p(v_1)=p(v_2)$, 
either $p(\sigma_{v_1})=p(\sigma_{v_2})$ or $p(\sigma_{v_1})=p(\sigma^{-1}_{v_2})$ holds. 

\bigskip 

Under this assumption, the pre-images of $v$ can be classified into two groups as follows. 
Assume that $V(G)=\{v_1, \ldots, v_h\}$ and 
$V(\tilde{G})=V_1\cup \cdots \cup V_h$ such that 
$v\in V_i$ if and only if $p(v)=v_i$ for $i\in \{1, \ldots, h\}$.  
Then, a {\em rotation-preserving coloring\/} 
$c \colon V(\tilde{G}) \to \{\mbox{black}, \mbox{white}\}$ is defined as follows. 
First, for each $i\in \{1, \ldots, h\}$, 
we fix a vertex $v\in V_i$, and color $v$ in black. 
For any other vertex $u$ of $V_i$, we define that 
$u$ is colored in black if and only if $p(\sigma_v)=p(\sigma_u)$; 
note that if $p(\sigma_v)=p(\sigma^{-1}_u)$, then $c(u)$ is white. 
Then, we decompose $V(\tilde{G})$ into $B\cup W$, where 
$B$ (resp. $W$) stands for the set of black (resp. white) colored 
vetices of $V(\tilde{G})$. 
For a rotation-preserving coloring $c$ of $\tilde{G}$, 
an edge $e=uv$ is {\em synchronous} if $c(u) = c(v)$, and {\em anti-synchronous} 
otherwise. Note further that this classification of edges 
also depends on $c$. 

\bigskip

\noindent
\textbf{Assumption 2:} 
For each edge $uv\in E(G)$, either every element of $p^{-1}(uv)$ is synchronous 
or every element of $p^{-1}(uv)$ is anti-synchronous.

\bigskip

As mentioned above, there are generally several rotation-preserving colorings on $\tilde{G}$. 
For a fixed $i \in \{1, \ldots, h\}$, exchanging the black and white colors 
assigned to $V_i \subset V(\tilde{G})$ switches the status of any edge 
between $V_i$ and $V_j$ $(i \neq j)$ from synchronous to anti-synchronous, and vice versa. 
However, the rotation compatibility of $(\tilde{G}, \sigma)$ remains unchanged. 
That is, this property does not depend on the choice of rotation-preserving coloring.

Negami showed the following theorem for rotation compatible covers.

\begin{thm}[Negami \cite{Negami2024main}]\label{thm:main_result}
A connected graph $G$ has a finite planar cover which is rotation compatible 
if and only if $G$ can be embedded on the projective plane.
\end{thm}

By Theorem~\ref{thm:main_result}, if every connected graph that admits a finite planar 
cover also has a finite rotation-compatible cover, then the Planar Cover Conjecture is solved. 
This approach, as noted above, differs from the conventional scheme,
which primarily focuses on showing that $K_{1,2,2,2}$ has no finite planar cover.
The proof of Theorem~\ref{thm:main_result} in \cite{Negami2024main}
relies on advanced algebraic and topological methods.
In this paper, we provide another simple proof based on combinatorial arguments,
which will be presented in the next section.
%

\section{Another simple proof of Theorem 2}\label{sect:Proof}

Let $G$ be a connected graph and let $p\colon \tilde{G} \to G$ be its $n$-fold cover.
Assume that $\tilde{G}$ is embedded on the plane and is rotation compatible.
Let $T$ be a spanning tree of $G$, which we call a {\em fixed spanning tree} in our argument.
Then $\tilde{G}$ contains $n$ disjoint trees $ T_1, T_2, \dots, T_n $ as subgraphs, each of which projects to $T$ by the covering map $p$; we put $\mathcal{T} = \{T_1, T_2, \dots, T_n\}$. 
Next, we consider the rotation-preserving coloring of $\tilde{G}$ 
with respect to the fixed tree $T$, denoted by $c_T $, such that each vertex of 
$T_1$ receives black by $c_T $. 
Then, for each $i \in \{2, 3, \dots, n\}$, 
vertices of $T_i$ are assigned the same color (either black or white) by Assumption~2; 
that is, every tree in $\mathcal{T}$ is monochromatic by $c_T $. 
Therefore, we have the decomposition 
$\mathcal{T} = \mathcal{T}_B \cup \mathcal{T}_W$, such that 
$\mathcal{T}_B$ (resp.\ $\mathcal{T}_W$) denotes the set of black colored 
(resp. white colored) trees by $c_T$.

For each $T_i \in \mathcal{T}$, we regard the induced subgraph
$\tilde{G}[V(T_i)]$ as a subembedding of $\tilde{G}$; that is, it 
is equipped with the embedding inherited from $\tilde{G}$. 
If every finite face of $\tilde{G}[V(T_i)]$ contains no vertex of $\tilde{G}$, then $T_i$ is \textit{innermost}.
Observe that $\tilde{G}$ admits at least one innermost tree in $\mathcal{T}$ since $T_1, T_2, \dots, T_n$ are disjoint.
We arbitrarily choose an innermost tree $T_k$ in $\mathcal{T}$ and call it the {\em core tree}. 
Without loss of generality, we may assume that $T_1$ is the core tree throughout. 
We define an \textit{enclosing curve} $\gamma$ of the core tree $T_1$ as a 
simple closed curve on the plane enclosing $\tilde{G}[V(T_1)]$ such that the disc bounded by 
$\gamma$ contains no vertices of $T_i$ ($i \neq 1$).
We can fix $\gamma$ such that for every edge intersecting $\gamma$, one endpoint belongs to $V(T_1)$ 
and the other belongs to $V(T_i)$ for some $i \neq 1$.
The edges of $\tilde{G}$ intersecting the enclosing curve $\gamma$ are 
{\em $\gamma$-crossing edges}, and the set of these edges is denoted by $E_\gamma$.

Let $e_1, e_2 \in E_\gamma$ be two edges mapped to the same edge $uv$ of $G$ 
by the covering map $p$.
We call the set $\{e_1, e_2\}$ a $\gamma$-{\em crossing pair} and denote it by $\gamma(u,v)$. 
Observe that for any edge $e \in E_\gamma$, there exists a unique edge $e' \in 
E_\gamma \setminus \{e\}$ such that $p(e)=p(e')$; that is, $\{e,e'\}$ is a 
$\gamma$-crossing pair. 
This implies that $|E_\gamma|$ is even. 
For a $\gamma$-crossing pair $\{e_1, e_2\}$, if $e_1$ is synchronous (resp. anti-synchronous), 
then $e_2$ is also synchronous (resp. anti-synchronous) by Assumption~2.
Thus, we use the terms ``synchronous'' and ``anti-synchronous'' for $\gamma$-crossing pairs.

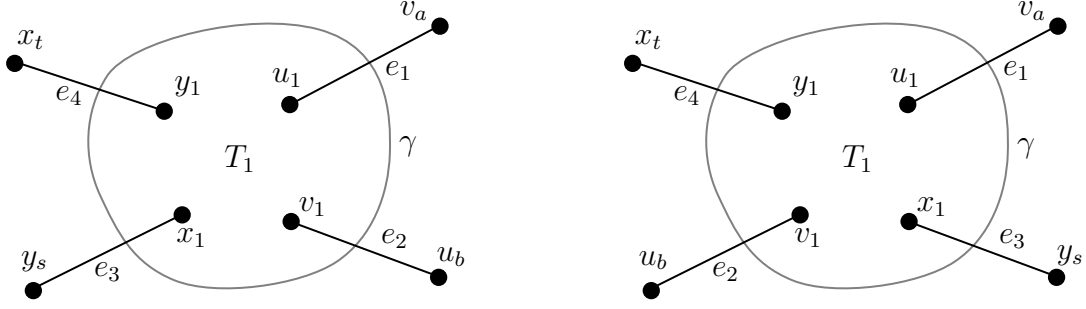
\begin{figure}[thb] 
\begin{center}

\tikzset{every picture/.style={line width=0.75pt}} 

\begin{tikzpicture}[x=0.75pt,y=0.75pt,yscale=-1,xscale=1]

\draw  [color={rgb, 255:red, 128; green, 128; blue, 128 }  ,draw opacity=1 ] (93.83,154.83) .. controls (105.47,136.97) and (165,122) .. (199.97,131.47) .. controls (234.94,140.94) and (236.33,174) .. (235.67,194) .. controls (235,214) and (227,236) .. (209,248) .. controls (191,260) and (148.67,265.83) .. (127.67,257.33) .. controls (106.67,248.83) and (99,233.33) .. (90.33,214.67) .. controls (81.67,196) and (82.2,172.7) .. (93.83,154.83) -- cycle ;
\draw [line width=0.75]    (130,225.02) -- (54.33,262.67) ;
\draw [line width=0.75]    (261,129.33) -- (183.61,170.23) ;
\draw [line width=0.75]    (261.67,256) -- (184.33,227.33) ;
\draw [line width=0.75]    (122.27,172.23) -- (44.33,146.67) ;
\draw  [fill={rgb, 255:red, 0; green, 0; blue, 0 }  ,fill opacity=1 ] (118.67,172.5) .. controls (118.67,170.38) and (120.38,168.67) .. (122.5,168.67) .. controls (124.62,168.67) and (126.33,170.38) .. (126.33,172.5) .. controls (126.33,174.62) and (124.62,176.33) .. (122.5,176.33) .. controls (120.38,176.33) and (118.67,174.62) .. (118.67,172.5) -- cycle ;
\draw  [fill={rgb, 255:red, 0; green, 0; blue, 0 }  ,fill opacity=1 ] (127.67,224.5) .. controls (127.67,222.38) and (129.38,220.67) .. (131.5,220.67) .. controls (133.62,220.67) and (135.33,222.38) .. (135.33,224.5) .. controls (135.33,226.62) and (133.62,228.33) .. (131.5,228.33) .. controls (129.38,228.33) and (127.67,226.62) .. (127.67,224.5) -- cycle ;
\draw  [fill={rgb, 255:red, 0; green, 0; blue, 0 }  ,fill opacity=1 ] (53,262.5) .. controls (53,260.38) and (54.72,258.67) .. (56.83,258.67) .. controls (58.95,258.67) and (60.67,260.38) .. (60.67,262.5) .. controls (60.67,264.62) and (58.95,266.33) .. (56.83,266.33) .. controls (54.72,266.33) and (53,264.62) .. (53,262.5) -- cycle ;
\draw  [fill={rgb, 255:red, 0; green, 0; blue, 0 }  ,fill opacity=1 ] (43.67,148.5) .. controls (43.67,146.38) and (45.38,144.67) .. (47.5,144.67) .. controls (49.62,144.67) and (51.33,146.38) .. (51.33,148.5) .. controls (51.33,150.62) and (49.62,152.33) .. (47.5,152.33) .. controls (45.38,152.33) and (43.67,150.62) .. (43.67,148.5) -- cycle ;
\draw  [fill={rgb, 255:red, 0; green, 0; blue, 0 }  ,fill opacity=1 ] (181.67,169.17) .. controls (181.67,167.05) and (183.38,165.33) .. (185.5,165.33) .. controls (187.62,165.33) and (189.33,167.05) .. (189.33,169.17) .. controls (189.33,171.28) and (187.62,173) .. (185.5,173) .. controls (183.38,173) and (181.67,171.28) .. (181.67,169.17) -- cycle ;
\draw  [fill={rgb, 255:red, 0; green, 0; blue, 0 }  ,fill opacity=1 ] (257,129.83) .. controls (257,127.72) and (258.72,126) .. (260.83,126) .. controls (262.95,126) and (264.67,127.72) .. (264.67,129.83) .. controls (264.67,131.95) and (262.95,133.67) .. (260.83,133.67) .. controls (258.72,133.67) and (257,131.95) .. (257,129.83) -- cycle ;
\draw  [fill={rgb, 255:red, 0; green, 0; blue, 0 }  ,fill opacity=1 ] (256.33,255.83) .. controls (256.33,253.72) and (258.05,252) .. (260.17,252) .. controls (262.28,252) and (264,253.72) .. (264,255.83) .. controls (264,257.95) and (262.28,259.67) .. (260.17,259.67) .. controls (258.05,259.67) and (256.33,257.95) .. (256.33,255.83) -- cycle ;
\draw  [fill={rgb, 255:red, 0; green, 0; blue, 0 }  ,fill opacity=1 ] (182.33,227.83) .. controls (182.33,225.72) and (184.05,224) .. (186.17,224) .. controls (188.28,224) and (190,225.72) .. (190,227.83) .. controls (190,229.95) and (188.28,231.67) .. (186.17,231.67) .. controls (184.05,231.67) and (182.33,229.95) .. (182.33,227.83) -- cycle ;
\draw  [color={rgb, 255:red, 128; green, 128; blue, 128 }  ,draw opacity=1 ] (403.83,154.83) .. controls (415.47,136.97) and (475,122) .. (509.97,131.47) .. controls (544.94,140.94) and (546.33,174) .. (545.67,194) .. controls (545,214) and (537,236) .. (519,248) .. controls (501,260) and (458.67,265.83) .. (437.67,257.33) .. controls (416.67,248.83) and (409,233.33) .. (400.33,214.67) .. controls (391.67,196) and (392.2,172.7) .. (403.83,154.83) -- cycle ;
\draw [line width=0.75]    (440,225.02) -- (364.33,262.67) ;
\draw [line width=0.75]    (571,129.33) -- (493.61,170.23) ;
\draw [line width=0.75]    (571.67,256) -- (494.33,227.33) ;
\draw [line width=0.75]    (432.27,172.23) -- (354.33,146.67) ;
\draw  [fill={rgb, 255:red, 0; green, 0; blue, 0 }  ,fill opacity=1 ] (428.67,172.5) .. controls (428.67,170.38) and (430.38,168.67) .. (432.5,168.67) .. controls (434.62,168.67) and (436.33,170.38) .. (436.33,172.5) .. controls (436.33,174.62) and (434.62,176.33) .. (432.5,176.33) .. controls (430.38,176.33) and (428.67,174.62) .. (428.67,172.5) -- cycle ;
\draw  [fill={rgb, 255:red, 0; green, 0; blue, 0 }  ,fill opacity=1 ] (437.67,224.5) .. controls (437.67,222.38) and (439.38,220.67) .. (441.5,220.67) .. controls (443.62,220.67) and (445.33,222.38) .. (445.33,224.5) .. controls (445.33,226.62) and (443.62,228.33) .. (441.5,228.33) .. controls (439.38,228.33) and (437.67,226.62) .. (437.67,224.5) -- cycle ;
\draw  [fill={rgb, 255:red, 0; green, 0; blue, 0 }  ,fill opacity=1 ] (363,262.5) .. controls (363,260.38) and (364.72,258.67) .. (366.83,258.67) .. controls (368.95,258.67) and (370.67,260.38) .. (370.67,262.5) .. controls (370.67,264.62) and (368.95,266.33) .. (366.83,266.33) .. controls (364.72,266.33) and (363,264.62) .. (363,262.5) -- cycle ;
\draw  [fill={rgb, 255:red, 0; green, 0; blue, 0 }  ,fill opacity=1 ] (353.67,148.5) .. controls (353.67,146.38) and (355.38,144.67) .. (357.5,144.67) .. controls (359.62,144.67) and (361.33,146.38) .. (361.33,148.5) .. controls (361.33,150.62) and (359.62,152.33) .. (357.5,152.33) .. controls (355.38,152.33) and (353.67,150.62) .. (353.67,148.5) -- cycle ;
\draw  [fill={rgb, 255:red, 0; green, 0; blue, 0 }  ,fill opacity=1 ] (491.67,169.17) .. controls (491.67,167.05) and (493.38,165.33) .. (495.5,165.33) .. controls (497.62,165.33) and (499.33,167.05) .. (499.33,169.17) .. controls (499.33,171.28) and (497.62,173) .. (495.5,173) .. controls (493.38,173) and (491.67,171.28) .. (491.67,169.17) -- cycle ;
\draw  [fill={rgb, 255:red, 0; green, 0; blue, 0 }  ,fill opacity=1 ] (567,129.83) .. controls (567,127.72) and (568.72,126) .. (570.83,126) .. controls (572.95,126) and (574.67,127.72) .. (574.67,129.83) .. controls (574.67,131.95) and (572.95,133.67) .. (570.83,133.67) .. controls (568.72,133.67) and (567,131.95) .. (567,129.83) -- cycle ;
\draw  [fill={rgb, 255:red, 0; green, 0; blue, 0 }  ,fill opacity=1 ] (566.33,255.83) .. controls (566.33,253.72) and (568.05,252) .. (570.17,252) .. controls (572.28,252) and (574,253.72) .. (574,255.83) .. controls (574,257.95) and (572.28,259.67) .. (570.17,259.67) .. controls (568.05,259.67) and (566.33,257.95) .. (566.33,255.83) -- cycle ;
\draw  [fill={rgb, 255:red, 0; green, 0; blue, 0 }  ,fill opacity=1 ] (492.33,227.83) .. controls (492.33,225.72) and (494.05,224) .. (496.17,224) .. controls (498.28,224) and (500,225.72) .. (500,227.83) .. controls (500,229.95) and (498.28,231.67) .. (496.17,231.67) .. controls (494.05,231.67) and (492.33,229.95) .. (492.33,227.83) -- cycle ;

\draw (238.89,183.77) node [anchor=north west][inner sep=0.75pt]  [font=\normalsize]  {$\gamma $};
\draw (151.17,188.07) node [anchor=north west][inner sep=0.75pt]    {$T_{1}$};
\draw (175.17,149.4) node [anchor=north west][inner sep=0.75pt]    {$u_{1}$};
\draw (188.5,214.07) node [anchor=north west][inner sep=0.75pt]    {$v_{1}$};
\draw (127.17,229.73) node [anchor=north west][inner sep=0.75pt]    {$x_{1}$};
\draw (126.17,153.4) node [anchor=north west][inner sep=0.75pt]    {$y_{1}$};
\draw (239.17,115.73) node [anchor=north west][inner sep=0.75pt]    {$v_{a}$};
\draw (257.83,239.4) node [anchor=north west][inner sep=0.75pt]    {$u_{b}$};
\draw (47.17,130.4) node [anchor=north west][inner sep=0.75pt]    {$x_{t}$};
\draw (49.17,241.07) node [anchor=north west][inner sep=0.75pt]    {$y_{s}$};
\draw (231.83,146.07) node [anchor=north west][inner sep=0.75pt]    {$e_{1}$};
\draw (229.83,230.4) node [anchor=north west][inner sep=0.75pt]    {$e_{2}$};
\draw (85.83,246.73) node [anchor=north west][inner sep=0.75pt]    {$e_{3}$};
\draw (66.5,158.73) node [anchor=north west][inner sep=0.75pt]    {$e_{4}$};
\draw (548.89,184.44) node [anchor=north west][inner sep=0.75pt]  [font=\normalsize]  {$\gamma $};
\draw (461.17,188.07) node [anchor=north west][inner sep=0.75pt]    {$T_{1}$};
\draw (485.17,149.4) node [anchor=north west][inner sep=0.75pt]    {$u_{1}$};
\draw (498.5,214.07) node [anchor=north west][inner sep=0.75pt]    {$x_{1}$};
\draw (437.17,229.73) node [anchor=north west][inner sep=0.75pt]    {$v_{1}$};
\draw (436.17,153.4) node [anchor=north west][inner sep=0.75pt]    {$y_{1}$};
\draw (549.17,115.73) node [anchor=north west][inner sep=0.75pt]    {$v_{a}$};
\draw (568.83,238.4) node [anchor=north west][inner sep=0.75pt]    {$y_{s}$};
\draw (357.17,130.4) node [anchor=north west][inner sep=0.75pt]    {$x_{t}$};
\draw (359.17,241.07) node [anchor=north west][inner sep=0.75pt]    {$u_{b}$};
\draw (541.83,146.07) node [anchor=north west][inner sep=0.75pt]    {$e_{1}$};
\draw (539.83,230.4) node [anchor=north west][inner sep=0.75pt]    {$e_{3}$};
\draw (395.83,246.73) node [anchor=north west][inner sep=0.75pt]    {$e_{2}$};
\draw (376.5,158.73) node [anchor=north west][inner sep=0.75pt]    {$e_{4}$};

\end{tikzpicture}
\caption{Consecutive $\gamma$-crossing pair and alternative $\gamma$-crossing pair.}
\label{fig:crossing_pair} 
\end{center} 
\end{figure}

Let $\gamma(u,v)= \{e_1, e_2\}$ and $\gamma(x,y)= \{e_3, e_4\}$ be two $\gamma$-crossing pairs. 
Then we define the following terms. 
\begin{itemize}
\item $\gamma(u,v)$ and $\gamma(x,y)$ are 
{\em consecutive} if edges $e_1, e_2, e_3, e_4$ (or $e_1, e_2, e_4, e_3$) cross $\gamma$ in this cyclic 
order (clockwise or anticlockwise); see the left-hand side of Figure~\ref{fig:crossing_pair}. 
\item $\gamma(u,v)$ and $\gamma(x,y)$ are {\em alternative} 
if edges $e_1, e_3, e_2, e_4$ cross $\gamma$ in this cyclic order (clockwise or anticlockwise); 
see the right-hand side of Figure~\ref{fig:crossing_pair}. 
\end{itemize}

Here, we show the key lemma for the proof of Theorem~\ref{thm:main_result}. 

\begin{lm}\label{lem:main}
For the core tree $T_1 \in \mathcal{T}$ of an $n$-fold rotation compatible planar cover 
$p \colon \tilde{G} \to G$, and for its enclosing curve $\gamma$, 
the following three properties hold:
\begin{itemize}
\item[{\rm (I)}]
Any two synchronous $\gamma$-crossing pairs are consecutive.
\item[{\rm (II)}]
Any two anti-synchronous $\gamma$-crossing pairs are alternative.
\item[{\rm (III)}]
A synchronous $\gamma$-crossing pair and an anti-synchronous $\gamma$-crossing pair are consecutive.
\end{itemize}
\end{lm}

\begin{Proof}
Without loss of generality, 
we may assume that 
$T_1 \in \mathcal{T}_B$.
In the following argument, for a vertex $v$ of $G$,
we denote the vertex of $p^{-1}(v) \cap V(T_i) \subset V(\tilde{G})$ by 
$v_i$ with the subscription 
for $i \in \{1, \dots, n\}$.
Let $\gamma(u, v) = \{ u_1v_a, v_1u_b \}$ and 
$\gamma(x, y) = \{ x_1y_s, y_1x_t \}$ be $\gamma$-crossing pairs, where $a, b, s, t \in \{2, \dots ,n\}$.
Let $P_{uv}$ and $P_{xy}$ be the unique paths in a spanning tree 
$T$ in $G$ joining $u$ and $v$, and $x$ and $y$, respectively. 
Then, 
$\mathcal{C}_{uv}=p^{-1}(P_{uv} \cup uv)$ and $\mathcal{C}_{xy}=p^{-1}(P_{xy} \cup xy)$ 
are unions of cycles of $\tilde{G}$. 
Let $C$ and $C'$ denote cycles in 
$\mathcal{C}_{uv}$ and $\mathcal{C}_{xy}$, 
respectively, that pass through $T_1$. 
Furthermore, we put $\mathcal{T}_{C} = \{T_i \in \mathcal{T} \colon V(C) \cap V(T_i) \neq \emptyset \}$, 
and let $D_C$ denote the $2$-cell region bounded by $C$. 
We assume that $C$ consists of $m$ edges of $p^{-1}(uv)$ and 
$m$ paths of $p^{-1}(P_{uv})$. 
Since $T_1 \in \mathcal{T}_{C}$, we relabel the trees in $\mathcal{T}$ as 
$\mathcal{T}_{C} = \{T_1, \dots ,T_m \}$ so that $T_1, \dots ,T_m$ intersect $C$ in this order. 
Then, $C$ can be expressed as
$$ C= P_{u_1v_1}v_1u_2P_{u_2v_2}v_2u_3P_{u_3v_3}v_3u_4 \cdots v_{m-1}u_mP_{u_mv_m}v_mu_1, $$
where $P_{u_iv_i} = p^{-1}(P_{uv}) \cap T_i$. 
We show the three cases (I)--(III) in turn. 

\smallskip
\smallskip

\noindent
Case (I): $\gamma(u, v)$ and $\gamma(x, y)$ are synchronous. 
Suppose for a contradiction that $\gamma(u, v)$ and $\gamma(x, y)$ are 
alternative. 
We may assume that $u_1v_a, x_1y_s, v_1u_b$, and $y_1x_t$ cross $\gamma$ in this clockwise order 
(see the right-hand side of Figure~\ref{fig:crossing_pair}). 
Since $T_1 \in \mathcal{T}_B$ and two edges in $\gamma(u, v)$ are synchronous, 
we have $\mathcal{T}_{C} \subset \mathcal{T}_B$. 
Since $\gamma(u, v)$ and $\gamma(x, y)$ are alternative, 
$C$ and $C'$ intersect at $T_l \in \mathcal{T}_{C}$ $(2 \leq l \leq m)$;
we may assume that $C'$ does not intersect $T_i$ for any $i \in \{2, \dots, l-1\}$ 
(see Figure~\ref{fig:Case(I)}). 
Then we can take a simple closed curve $\gamma'$ surrounding $T_l$ such that 
for every edge intersecting $\gamma'$, one endpoint belongs 
to $V(T_l)$ and the other belongs to $V(T_j)$ for some $j \neq l$. 
Under the condition, 
$u_lv_{l-1}, x_ly_r, v_lu_{l+1}$, and $y_lx_q$ 
$(q, r \in \{1, \ldots, n\}\setminus \{l\})$ 
cross $\gamma'$ in this clockwise order. 
This implies that $T_l \in \mathcal{T}_W$, which contradicts $\mathcal{T}_{C} \subset \mathcal{T}_B$.
Therefore, we got our desired conclusion. 

\begin{figure}[htb]
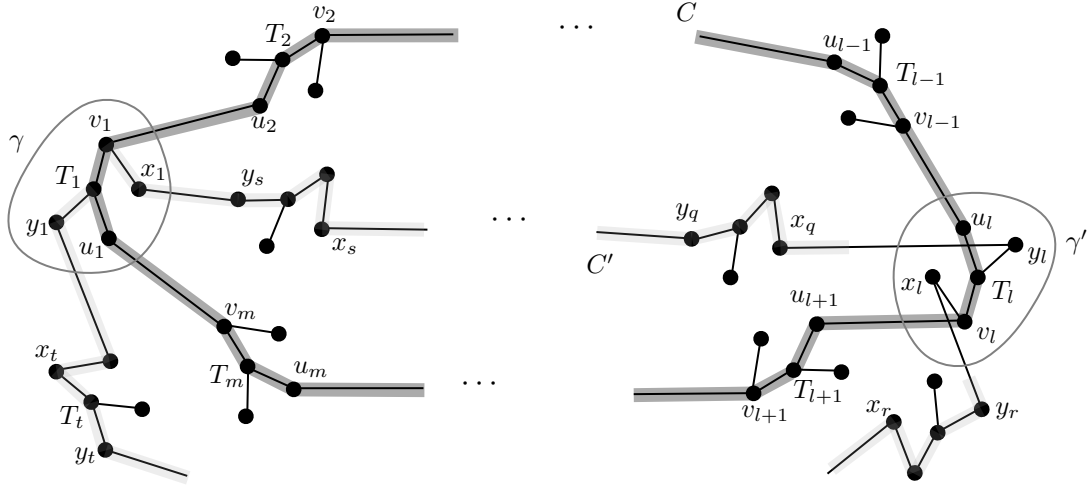
 
\begin{center}

\tikzset{every picture/.style={line width=0.75pt}} 


\caption{Configulation of Case (I).}
\label{fig:Case(I)} 
\end{center} 
\end{figure}

\smallskip
\smallskip

\noindent
Case (II): $\gamma(u, v)$ and $\gamma(x, y)$ are anti-synchronous. 
Suppose for a contradiction that $\gamma(u, v)$ and $\gamma(x, y)$ 
are consecutive. 
We may assume that $u_1v_a, v_1u_b, x_1y_s$, and $y_1x_t$ cross $\gamma$ in this clockwise order 
(see the left-hand side of Figure~\ref{fig:crossing_pair}). 
Since $\gamma(u, v)$ is anti-synchronous, $m$, which is the number of trees 
in $\mathcal{T}_C$, is even. 
Recall that $T_1 \in \mathcal{T}_B$, and hence for $i \in \{1, \ldots, m\}$, $T_i \in \mathcal{T}_W$ 
if and only if $i$ is even. 
Without loss of generality, 
we may assume that $x_1y_s$ and $y_1x_t$, which are of $\gamma(x, y)$, lie in $D_C$ (see Figure~\ref{fig:Case(II)}). 
Then for $i \in \{1, \ldots, m\}$, two edges of $p^{-1}(xy)$ which are incident to $T_i$ lie in 
$D_C$ if and only if $i$ is odd. 
We further divide our argument into two small cases (II--i) and (II--ii) below. 

\begin{figure}[htb]
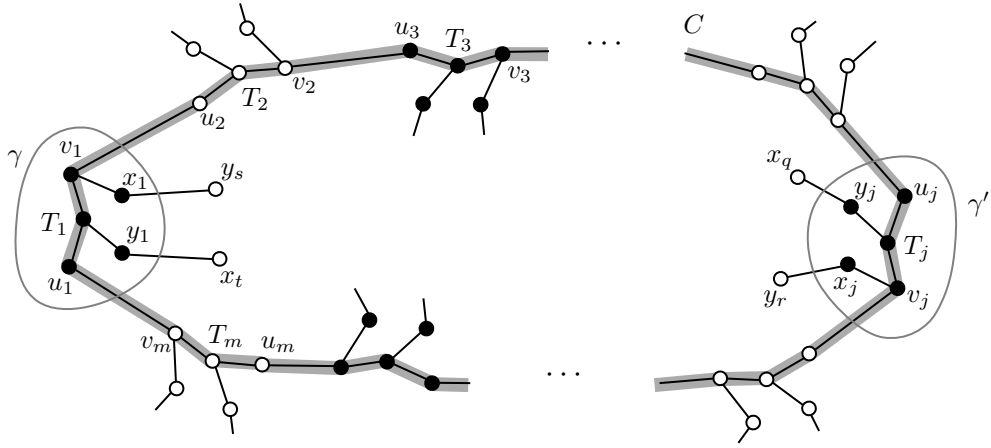
 
\begin{center}

\tikzset{every picture/.style={line width=0.75pt}} 


\caption{Configulation of Case (II).}
\label{fig:Case(II)} 
\end{center} 
\end{figure}

\smallskip

\noindent
(II--i)\quad
Consider the case where $D_C$ 
contains no tree of $\mathcal{T} \setminus \mathcal{T}_C$ in its interior.
Under our assumption, there exists $j \in \{2, \dots , n\}$ such that $x_1y_j$ is inside the region.
By the argument above, $j$ must be odd. 
However it contradicts that $\gamma(x, y)$ is anti-synchronous since $T_1,T_j \in \mathcal{T}_B$.

\smallskip

\noindent
(II--ii)\quad
Consider the case when $D_C$ 
contains some tree of $\mathcal{T} \setminus \mathcal{T}_C$ in its interior.
In this case, $D_C$ contains some cycles of $\mathcal{C}_{uv}$. 
We carry out the same argument as (II--i) for the innermost one, and obtain the contradiction.

\smallskip
\smallskip

\noindent
Case (III): one of $\gamma(u,v)$ and $\gamma(x,y)$ is synchronous and the other is anti-synchronous. 
The same argument as in the case (I) works.
Therefore, a synchronous $\gamma$-crossing pair and an anti-synchronous $\gamma$-crossing 
pair are consecutive.
\end{Proof}

Now, we prove our main result. 

\bigskip

\begin{Proofof}{Theorem \ref{thm:main_result}}
The sufficiency is clear, and hence we prove the necessity.
Let $\tilde{G}$ be a rotation compatible planar 
cover of $G$ with a covering map 
$p\colon \tilde{G} \to G$. 
Let $T_1 \in \mathcal{T}$ be the core tree in $\tilde{G}$ associated with 
a spanning tree $T$ in $G$, and let $\gamma$ be an enclosing curve of $\tilde{G}[V(T_1)]$.
Then we can decompose the set of $\gamma$-crossing edges $E_\gamma$ 
into the set of synchronous edges $S_\gamma$ and the set of anti-synchronous edges $A_\gamma$, 
that is, $E_\gamma = S_\gamma \cup A_\gamma$ and $S_\gamma \cap A_\gamma = \emptyset$. 
We can fix $\gamma$ so that for each edge $e$ in $E_\gamma$, $\gamma$ crosses $e$ 
at its {\em midpoint}, which is taken as an interior point of $e$. 

Let $D$ be the closed disk bounded by $\gamma$. 
Then, $\tilde{G}[V(T_1)]$ is embedded on $D$ without any edge crossing, 
and $D$ contains {\em half-edges} of edges in $E_\gamma$, whose 
endpoints on $\gamma$ are midpoints of those (full-)edges; 
we call those endpoints of half-edges that are not vertices of $\tilde{G}$ 
are {\em fake vertices}. 
For a $\gamma$-crossing pair $\gamma(u,v)$, we denote the 
two half-edges of $\gamma(u,v)$ in $D$ by $u_1f_{uv}$ and 
$v_1f_{vu}$, where $f_{xy}$ represents a fake vertex on an edge $x_1y_i$ 
whose 
one end vertex $y_i$ outside $D$ for some 
$i \in \{2, \ldots, n\}$ with $p(x_1)=x\in V(G)$ and $p(y_i)=y\in V(G)$ 
(see the left-hand side of Figure~\ref{fig:proof_disc}). 
Furthermore, we denote by $\hat{A}_\gamma$ and $\hat{S}_\gamma$ 
the set of the half-edges of $A_\gamma$ and 
$S_\gamma$ in $D$, respectively.

Now, consider the closed unit disk $D_0$ 
with its boundary partitioned into $|A_\gamma|$ arcs of equal length. 
Let $h_0, h_1, \ldots, h_{|A_\gamma|-1}$ be the endpoints of these arcs, 
labeled in cyclic order along the boundary. 
We can take a homeomorphism $\phi \colon D \to D_0$ such that 
all the fake vertices of $\hat{A}_\gamma$ correspond 
to $h_0, h_1, \ldots, h_{|A_\gamma|-1}$. 
Under the condition, for any anti-synchronous 
$\gamma$-crossing pair $\gamma(u,v)$, 
two fake vertices $f_{uv}$ and $f_{vu}$ on edges of $\gamma(u,v)$ 
placed at antipodal pair of points on $\phi(\gamma)$, that is, $h_i$ and $h_{i+|\frac{|A_\gamma|}{2}|}$
with subscripts taken modulo $|A_\gamma|$, by 
Lemma~\ref{lem:main} (II); 
we can easily show the fact by induction on $|A_\gamma|$. 
See the middle drawing of Figure~\ref{fig:proof_disc}. 
In the middle and right-hand drawings of this figure, 
we adopt the notation that
$a_i$ and $a_i'$ (resp. $s_i$ and $s_i'$) ($i \geq 1$) denote corresponding vertices of 
$\hat{A}_\gamma$ (resp. $\hat{S}_\gamma$) inside $D_0$.

\begin{figure}[tb]
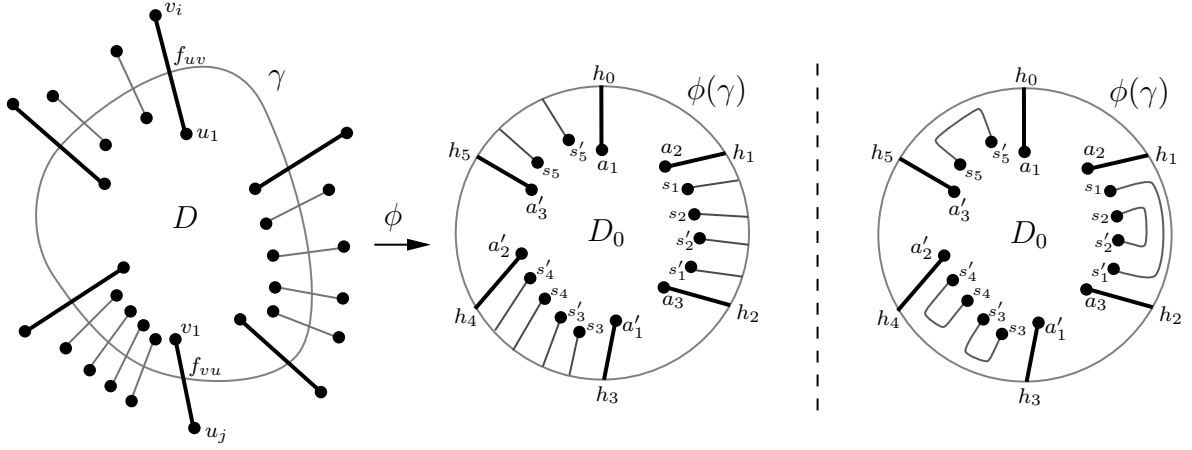
 
\begin{center}

\tikzset{every picture/.style={line width=0.75pt}} 


\caption{Edge arrangement in the disc.}
\label{fig:proof_disc} 
\end{center} 
\end{figure}

By Lemma~\ref{lem:main} (I), 
for any synchronous 
$\gamma$-crossing pair $\gamma(u,v)$, 
two fake vertices $f_{uv}$ and $f_{vu}$ on edges of $\gamma(u,v)$ 
placed on interior of some arc $h_ih_{i+1}$ on $\phi(\gamma)$. 
Furthermore, assume that there exists another 
synchronous $\gamma$-crossing pair $\gamma(x,y)$ such that 
its fake vertices $f_{xy}$ and $f_{yx}$ lie on the same arc $h_ih_{i+1}$. 
Then, the pairs $\{f_{uv}, f_{vu}\}$ and $\{f_{xy}, f_{yx}\}$ 
do not interleave on the arc by 
Lemma~\ref{lem:main} (III). 
This enables us to move fake vertices of $u_1f_{uv}, v_1f_{vu}, x_1f_{xy}$, 
and $y_1f_{yx}$ so as to identify $f_{uv}$ with $f_{vu}$ and $f_{xy}$ with $f_{yx}$, 
thereby obtaining two edges $u_1v_1$ and $x_1y_1$ in $D_0$ without any crossings.
This argument works for the case when $h_ih_{i+1}$ contains at least six fake vertices, 
also by Lemma~\ref{lem:main} (III); if necessary, use induction by removing the innermost pair. 
See the right-hand side of Figure~\ref{fig:proof_disc}. 

By argument above, by identifying each antipodal pair of points of the boundary of $D_0$, 
we obtain the projective plane in which the graph isomorphic to $G$ is embedded. 
\end{Proofof}

In this paper, we provide another simple proof of the main result in \cite{Negami2024main}, 
based purely on combinatorial arguments. 
However, when approaching the planar cover conjecture from this direction, 
the deep arguments in \cite{Negami2024main}, concerning regular covers 
and the construction of a (possibly punctured) closed surface as a quotient of 
a surface of larger genus---known as a branched cover of a surface---are essential. 
Together with our new proof, it is desirable to develop a deeper understanding 
of the structure of graph coverings. 

\section*{Acknowledgement}
This work is supported by JSPS KAKENHI Grant Number JP23K03196 
and JST SPRING, Grant Number JPMJSP2121.

\end{document}